\documentclass[11pt]{article}
\usepackage{graphicx}
\usepackage{amssymb}
\usepackage{epstopdf}
\def\qedbox{\hbox{$\rlap{$\sqcap$}\sqcup$}}
\def\qed{\nobreak\hfill\penalty250 \hbox{}\nobreak\hfill\qedbox}

\newtheorem{theorem}{Theorem}

\newtheorem{lemma}[theorem]{Lemma}

\title{Locally Isotropic Pseudo-Riemannian Manifolds}
\author{Iva Stavrov}
\begin{document}
\maketitle
\begin{abstract}
Locally isotropic pseudo-Riemannian manifolds are known to be locally symmetric; this result is due to Wolf (\cite{Wolf}). In the Riemannian setting one proof, due to Szab\'o, uses spectral properties of the so-called {\it Szab\'o operator}. In this paper we extend Szab\'o's method to the pseudo-Riemannian setting, obtaining results comparable to those of Wolf. Most results concerning the Szab\'o operator in the pseudo-Riemannian setting are presented in \cite{GilkeyReniJa}, primarily using methods of algebraic topology. This paper exploits the polynomial nature of the Szab\'o operator along with its behavior over the nullcone, both of which have received little attention this far.
\end{abstract}

\section{Introduction} 

\hskip .2in Let $M$ be a smooth manifold and let $(.,.)$ be a non-degenerate metric on the tangent bundle $TM$. Let $\nabla$ be the Levi-Civita connection on $TM$ and let $R$ be the corresponding Riemann curvature tensor. The covariant derivative of the curvature tensor $\nabla R$ is a section of the vector bundle $\otimes^5T^*M$ satisfying the following symmetries:
\begin{eqnarray}
&\nabla R(x,y,z,w;v)=\nabla R(z,w,x,y;v)= -\nabla R(y,x,z,w;v)\\
&\ \ \nabla R(x,y,z,w;v)+\nabla R(y,z,x,w;v)+\nabla R(z,x,y,w;v)=0\\
&\ \ \nabla R(x,y,z,w;v)+\nabla R(x,y,w,v;z)+\nabla R(x,y,v,z;w)=0.
\end{eqnarray}

Let $v$ be a tangent vector in the tangent space $T_PM$ at a point $P\in M$. The {\it Szab\'o operator} corresponding to $v$ is the operator $\mathbf{S}(v):T_PM\to T_PM$ defined by:
\begin{equation}
(\mathbf{S}(v)x,y)=\nabla R(x,v,v,y;v)\mathrm{\ for\ all\ } x,y\in T_PM.
\end{equation}
  
We are interested in the Szab\'o operator because it can be used in order to show that a certain manifold is locally symmetric using the following Theorem. 
\begin{theorem}
\label{thm2}
The Szab\'o operator $\mathbf{S}$ of a pseudo-Riemannian manifold $M$ vanishes identically if and only if the covariant derivative of the Riemann curvature tensor $\nabla R$ vanishes identically.
\end{theorem}

A purely algebraic proof  of this Theorem which uses nothing more than the curvature symmetries and polarization is given in Lemma 3.8.1 of \cite{elbook}. For more geometric ways to prove Theorem \ref{thm2} the reader is referred to 
\cite{bese}.

The Szab\'o operator at a point $P\in M$ can be viewed as a map 
$\mathbf{S}:T_PM\to \mathrm{Hom}(T_PM, T_PM)$.
The following properties of the map $\mathbf{S}$ follow easily from the symmetries stated in equations (1)-(3).

\begin {theorem}
\label{thm1}
Adopt the notation established above. Let $v\in T_PM$ be arbitrary, let $\mathcal{T}$ be a local isometry of $M$ which fixes $P$, and let $T:T_PM\to T_PM$ be the differential of $\mathcal{T}$ at $P$. We have:
\begin{enumerate}
  \item $\mathbf{S}(v)$ is self-adjoint, i.e. $(\mathbf{S}(v)x,y)=(x,\mathbf{S}(v)y)$ for all $x,y\in T_PM$; 
  \item $\mathbf{S}(v)v=0$;
  \item $\mathbf{S}(-v)=-\mathbf{S}(v)$;
  \item $\mathbf{S}(Tv)=T\circ\mathbf{S}(v)\circ T^{-1}$;
  \item If $\{e_1,e_2,\ldots,e_m\}$ is any basis for $T_PM$, then $\mathbf{S}(v_1e_1+v_2e_2+\ldots+v_me_m)$ is a homogeneous polynomial of degree $3$ in variables $v_1, v_2,\ldots,v_m$ with coefficients in $\mathrm{Hom}(T_PM,T_PM)$. In particular, the map $\mathbf{S}$ is continuous.
\end{enumerate}
\end {theorem}

We would like to point out that since we are in the pseudo-Riemannian setting, operators $\mathbf{S}(v)$ need not be diagonalizable. For details on this and related issues one faces in pseudo-Riemannian geometry the reader is referred to \cite{elbook}. 

Theorem \ref{thm1} motivates us to study the class $\mathcal{P}_P$ of maps $S:T_PM\to \mathrm{Hom}(T_PM,T_PM)$ which satisfy properties 1-5 of Theorem \ref{thm1}. We shall use 
$\mathcal{P}_{n,P}$ to denote the set of maps $S\in \mathcal{P}_P$ such that for any basis $\{e_1,e_2,\ldots,e_m\}$ of $T_PM$ we have that $S(v_1e_1+v_2e_2+\ldots+v_me_m)$ is a homogeneous polynomial of degree $2n+1$ in variables $v_1, v_2,\ldots,v_m$. When the basepoint $P$ is clear from the context we simply write $\mathcal{P}$ and $\mathcal{P}_n$. It can easily be checked that  classes $\mathcal{P}$ and $\mathcal{P}_{n}$ form real vector spaces for all non-negative integers $n$.  Note also  that the space $\mathcal{P}$ is closed under taking odd powers.

Theorem \ref{thm1} implies that the Szab\'o operator $\mathbf{S}$ is an element of $\mathcal{P}_1$.
In this paper we study the space $\mathcal{P}_1$, and in particular the Szab\'o operator, for a special class of manifolds. Following Wolf \cite{Wolf},
we say that  a pseudo-Riemannian manifold $M$ is {\it locally isotropic} if for any point $P$ and any two nonzero tangent vectors $x$ and $y$ at $P$ with $(x,x)=(y,y)$, there is a local isometry of $M$ fixing $P$, which sends $x$ to $y$. Wolf showed that locally isotropic manifolds are necessarily locally symmetric; see Theorem 12.3.1 of  \cite{Wolf}. In this paper we prove the following:
\begin{theorem}
\label{main}
Let $M$ be a locally isotropic pseudo-Riemannian manifold of signature $(p,q)$.
\begin{enumerate}
\item If $p=q$ and if $p\ne 2,4,8$, then $M$ is locally symmetric.
\item If $p\ne q$ with $\max\{p,q\}\ge 11$, then $M$ is locally symmetric.
\end{enumerate} 
\end{theorem}

In fact, we will show that if manifold $M$ is locally isotropic of signature $(p,q)$ satisfying one of the conditions of Theorem \ref{main}, then $\mathcal{P}_{1,P}$ is trivial for all $P\in M$. In particular, the Szab\'o operator $\mathbf{S}$ vanishes identically and we may conclude that manifold $M$ is locally symmetric, by Theorem \ref{thm2}. This method originates from the work of Szab\'o who in \cite{Szabo} used the spectral properties of the operator $\mathbf{S}$ in the case of a $2$-point homogeneous Riemannian manifold $M$ in order to show that $M$ is locally symmetric. 

While Szab\'o's method in the Riemannian setting relies heavily upon algebraic topology, it requires only elementary linear algebra to show that Lorentzian locally isotropic manifolds are locally symmetric. In fact,  Lorentizan locally isotropic manifolds have constant sectional curvature (see \cite{GilkeyStavrov}). The result  was obtained by studying the Jacobi operator; the crucial step in the proof being the observation that the Jacobi operator is nilpotent over the nullcone bundle. 
This suggests using the behavior of the Szab\'o operator over the nullcone bundle, along with methods of algebraic topology already common in the subject (see, for example, \cite{chi}, \cite{elbook}, \cite{GLS}, \cite{Za}). Most of the previous results concerning the Szab\'o operator are presented in \cite{GilkeyReniJa}, primarily using algebraic topology. The polynomial nature of the operators defined by the Riemann curvature tensor, although used earlier when studying the Jacobi operator (see \cite{Nikolaj}), has received little attention thus far. In this paper we exploit this polynomial nature in combination with nilpotency over the nullcone and the algebraic topology approach.

Here is the overview of the paper. In Section 2 we develop the necessary technical material from algebraic topology. In particular, we prove results about vector bundles over real projective spaces induced by elements of $\mathcal{P}_n$. We then prove various results regarding such vector bundles. In Section 3 we introduce vector valued polynomial maps and prove a number of lemmas which allow us to take advantage of the polynomial nature of the Szab\'o operator. In Section 4 we prove Wolf's Theorem for manifolds of signature $(p,p)$, where $p\ne 2,4,8$. It should be pointed out that the basis of our proof is the nilpotency result for elements of $\mathcal{P}$; this is an extension of the earlier nilpotency result of Gilkey-Stavrov \cite{GilkeyStavrov}. Wolf's Theorem in the general case is proved in Section 5, beginning with work of \cite{GilkeyReniJa}. We then use this, along with the polynomial nature of the Szab\'o operator, to construct a vector bundle of the type discussed in Section 2. The proof follows from the results of that section. 

\section{Background in Algebraic Topology}

\hskip .2in As pointed out in the Introduction, elements of $\mathcal{P}_n$ give rise to vector bundles over real projective spaces. It is for this reason that most of our computations involve Stiefel-Whitney classes of vector bundles. We now state the axioms of Stiefel-Whitney classes; for further details the reader is referred to \cite{MilnorStasheff}.

\begin{theorem} 
\label{SW}
To every vector bundle $V$ over a space $X$ we can associate 
an element $$w(V)=1+w_1(V)+\ldots+w_k(V)+\ldots \in H^*(X;\mathbb{Z}_2),$$ with homogeneous components $w_i(V)\in H^i(X;\mathbb{Z}_2)$, such that:
\begin{enumerate}
\item We have $w_k(V)=0$, when $k> \mathrm{rank}(V)$; 
\item For two vector bundles $V_1$ and $V_2$ we have
$w(V_1\oplus V_2)=w(V_1)w(V_2)$;
\item If $f:X\to Y$, then $w(f^*(V))=f^*\big{(}w(V)\big{)}$;
\item The element $w_1(\gamma_1)$ generates
$H^1(\mathbb{R}P^n;\mathbb{Z}_2)\cong\mathbb{Z}_2$.
\end{enumerate}
\end{theorem}

It follows from property 4 of Stiefel-Whitney classes that the generator of the truncated polynomial ring 
$H^*(\mathbb{R}P^{n};\mathbb{Z}_2)$ is the first Stiefel-Whitney class of the canonical line bundle, 
$x:=w_1(\gamma_1)$. In other words, 
$$H^*(\mathbb{R}P^{n};\mathbb{Z}_2)\cong \mathbb{Z}_2[x]/(x^{n+1}).$$
\bigbreak
We would now like to review the result of Adams \cite{A} regarding the K-theory of real projective spaces.

\begin {theorem}
\label{Adams}
Let  $\phi(n)$ denote the number of integers $s$ which satisfy 
$$1\le s\le n {\ \ and\ \ } s\equiv 0,\ 1,\ 2,\ 4\ \mathrm{mod}\ 8.$$
The stable equivalence class $\{\gamma_1\}\in \widetilde{KO}(\mathbb{R}P^n)$
generates
$\widetilde{KO}(\mathbb{R}P^n)$, and is of order $2^{\phi(n)}$.
\end{theorem}

Inspection shows that  $\frac{n-1}{2}\le\phi(n)\le\frac{n+2}{2}$. Given $n$, let $j$ be the unique integer satisfying $2^j\le n< 2^{j+1}$. We have that $j$ is roughly equal to $\log_2 n$. Therefore, for large values of $n$ we have $\phi(n)\ge j+3$. This technical inequality plays a significant role in the topological part of the proof of Wolf's Theorem. More careful treatment of the inequality is given in the following Lemma. 

\begin{lemma}
\label{phi}
Adopt the notation established above. 
\begin{enumerate}
\item If $n\ge 10$, then $\phi(n)\ge j+3$;
\item If $n\ne 1,3,7,$ then $2^{\phi(n)}>n+1$.
\end{enumerate}
\end{lemma}

\medbreak\noindent{\bf Proof.} The function $f(x)=\frac{x-1}{2}-\log_2 x$ is increasing for $x\ge 3$; this can easily be seen from the first derivative $f^\prime(x)$. Since $f(13)>2$, we have $\frac{n-1}{2}>\log_2 n+2$ for all $n\ge 13$. Inequalities $\phi(n)\ge \frac{n-1}{2}$ and $\log_2 n\ge j$ now imply 
$$\phi(n)>j+2, \mathrm{\ i.e.\ } \phi(n)\ge j+3 \mathrm{\ for\ all\ } n\ge 13.$$
Direct verification shows $\phi(n)\ge j+3$ for all $n\ge 10$. Combining the inequality $\phi(n)\ge j+3$ with $2^{j+1}>n$ gives us $2^{\phi(n)}> n+1$ for all $n\ge 10$. One now checks that $2^{\phi(n)}>n+1$ holds for all $n\ne 1,3,7$.\qed
\bigbreak 
Let $M$ be a locally isotropic pseudo-Riemannian manifold of signature $(p,q)$. Let $S\in \mathcal{P}_{P}$.
If $x$ and $y$ are two unit spacelike  (two unit timelike or two non-zero null) vectors at $P$, then there is a map $T:T_PM\to T_PM$ such that 
$$S(y)=T\circ S(x)\circ T^{-1}.$$
Thus, the rank of the operator $S(x)$ is independent of the choice of unit spacelike (resp. timelike or non-zero null) vector $x$.

In general, let $U$ be a finite dimensional vector space. A continuous map $S:X\to \mathrm{Hom}(U,U)$ such that the rank of $S(x)$ is independent of the choice of $x\in X$, gives rise to a vector bundle $E$ over $X$; the fibers of $E$ are determined by 
$$E\big{|}_x=\mathrm{ Im}\big{(}S(x)\big{)}.$$
An example of such a continuous map is $$S=\mathbf{S}\circ \iota:S^{q-1}\to \mathrm{Hom}(T_PM,T_PM),$$ where $\iota:S^{q-1}\to T_PM$ is the natural inclusion of the unit sphere $S^{q-1}$ into a maximal positive definite subspace of $T_PM$. In addition to knowing that the rank of $S(x)$ does not change with $x$, we also know that:
\begin{itemize}
\item $S(-x)=-S(x)$ for all $x$;
\item the operators $S(x)$ are self-adjoint for all $x$.
\end{itemize}
As a consequence, we are able to say more about the induced vector bundle $E$.

\begin{lemma}
\label{oddmaps}
Let $U$ be a finite dimensional vector space equipped with a non-degenerate inner product $(.,.)$. Consider a continuous map $S:S^{n}\to \mathrm{Hom}(U, U)$  such that: 
\begin{enumerate}
\item The rank of the operator $S(x)$ is the same for all $x\in S^{n}$;
\item We have $S(-x)=-S(x)$ for all $x\in S^{n}$;
\item The operators $S(x)$ are self-adjoint for all $x\in S^n$.
\end{enumerate}
Let $\pi:S^n\to \mathbb{R}P^n$ be the natural projection and let $E$ denote the vector bundle over $S^n$ with fibers $E\big{|}_x=\mathrm{ Im}\big{(}S(x)\big{)}.$ Then, there exists a vector bundle $\mathrm{Im}(S)$ over $\mathbb{R}P^n$ such that:
\begin{enumerate}
\item $\pi^*\mathrm{Im}(S)\cong E$;
\item $\mathrm{Im}(S)$ is a sub-bundle of the trivial bundle $\mathbb{R}P^n\times U$;
\item $\mathrm{Im}(S)\cong \mathrm{Im}(S)\otimes \gamma_1$;
\item If the rank of $\mathrm{Im}(S)$ is non-zero, then $\mathrm{Im}(S)$ allows a nowhere vanishing global section.
\end{enumerate}
\end{lemma}

\medbreak\noindent{\bf Proof.}
We see from  ${\rm Im}\big{(}S(x)\big{)}={\rm Im}\big{(}S(-x)\big{)}$ that the vector bundle $E$ descends to a vector bundle over $\mathbb{R}P^{n}$; we will denote this vector bundle by ${\rm Im}(S)$. Since $E$ is a sub-bundle of $S^n\times U$, we have that $\mathrm{Im}(S)$ is a sub-bundle of $\mathbb{R}P^n\times U$; assertions 1 and 2 now follow.

To prove the remaining two assertions note that we may, for the purposes of studying the vector bundle ${\rm Im}(S)$, assume the non-degenerate inner product $(.,.)$ on $U$ is actually positive definite. Indeed, let $U=U_-\oplus U_+$ be a decomposition of $U$ into a
direct sum of a maximal positive definite subspace $U_+$ and its orthogonal
complement $U_-$. Let
$\varrho_+:U\to U_+$ and $\varrho_-:U\to U_-$ denote the corresponding
orthogonal projections. Consider the linear map $\Phi:U\to U$ defined by
$\Phi v:=\varrho_+ v-\varrho_- v$ and the positive definite inner product $g$
defined by:
\begin{equation}
\label{eqn3}
\label{eqn4}
g(v,w):=(v,\Phi w)=(\Phi v,w).
\end{equation}
The correspondence  $T\mapsto T\circ \Phi$
is a bijection between the set of those operators on $U$ which are self-adjoint with respect to $(.,.)$ and  those which are self-adjoint with respect to $g$. Consequently, if $S:S^n\to \mathrm{Hom}(U,U)$ satisfies the conditions of the Lemma, so does $$\tilde{S}:S^n\to \mathrm{Hom}\big{(}(U,g),(U,g)\big{)} \mathrm{\ \ defined\ by\ \ } \tilde{S}(x):=S(x)\circ\Phi.$$  Since $\mathrm{Im}\big{(}\tilde{S}(x)\big{)}=\mathrm{Im}\big{(}S(x)\big{)},$ replacing $S$ by $\tilde{S}$ does not change the induced vector bundle. Therefore, in what follows we assume the inner product on $U$ is positive definite. 

The maps $S(x)$ are self-adjoint and so $\mathrm{Ker}\big{(}S(x)\big{)}\cap \mathrm{Im}\big{(}S(x)\big{)}=\{0\}$. As a consequence, $S(x)$ is an automorphism of $\mathrm{Im}\big{(}S(x)\big{)}$ for all $x\in S^n$. Thus we have a vector bundle isomorphism $\mathcal{S}:E\to E$. However, since $S(-x)=-S(x)$, the isomorphism $\mathcal{S}$ does not descend to an isomorphism of $\mathrm{Im}(S)$. Rather, it gives rise to a vector bundle isomorphism $$\mathrm{Im}\big{(}S(x)\big{)}\cong\mathrm{Im}\big{(}S(x)\big{)}\otimes\gamma_1.$$

We now prove the last assertion of the Lemma. The span of the eigenvectors corresponding to positive (resp. negative) eigenvalues of $S$ gives rise to a vector bundle $E_+$ (resp. $E_-$) over $S^n$; see Lemma 4.2.6 of \cite{elbook} for details. Note that if $a:S^n\to S^n$ is the antipodal map, then $a^*E_+\cong E_-$; this is due to the identity $S(-x)=-S(x)$. 
Let $A\in S^{n}$. Since $S^{n}-\{A\}$ is contractible, the vector bundle $E_+$ is trivial over $S^{n}-\{A\}$. Thus, there exists a nowhere vanishing section $e$ of $E_+$, defined over $S^{n}-\{A\}$. Multiplying by a smooth function on $S^{n}$ which vanishes only at $A$, we may assume that
the section $e$ is defined over the whole $S^{n}$ and vanishing only at  $A$. Consider the section $a^*e$ of $E_-$ corresponding to $e$ under the isomorphism $a^*E_+\cong E_-$; this section vanishes only at $a(A)$. Since $E_+\big{|}_x$ is orthogonal to $E_-\big{|}_x$ for all $x$, we see that the section $e+a^*e$ is nowhere vanishing over the entire sphere $S^n$. Moreover, 
$$[e+a^*e](-x)=e(-x)+[a^*e](-x)=[a^*e](x)+e(x)=[e+a^*e](x)$$
and so the section $e+a^*e$ of $E$ descends to a nowhere vanishing section of $\mathrm{Im}(S)$.\qed
\medbreak
In our work we often encounter isomorphisms as in property 3 of the previous Lemma. For this reason we often have to deal with vector bundles satisfying $V\cong V\otimes\gamma_1$ or with vector bundles of the form $V\oplus(V\otimes\gamma_1)$. In the following Lemma we study stable equivalence classes of such vector bundles.

\begin{lemma}
\label{lema5}
Let $V$ be a vector bundle over $\mathbb{R}P^{n}$ and let $\phi(n)$ be as in Theorem \ref{Adams}.
\begin{enumerate}
\item If $r=\mathrm{rank}(V)$, then $\{V\}+\{V\otimes \gamma_1\}=r\{\gamma_1\}$ in $\widetilde{KO}(\mathbb{R}P^{n})$;
\item If in addition $V\cong V\otimes \gamma_1$, then $\{V\}=a\{\gamma_1\}$ with $2a\equiv r \mathrm{\ mod\ }2^{\phi(n)}$.
\end{enumerate}
\end{lemma}

\medbreak\noindent{\bf Proof.}
Since $\widetilde{KO}(\mathbb{R}P^{n})$ is generated by the stable equivalence class $\{\gamma_1\}$, there exists an integer $a$ such that $\{V\}=a\{\gamma_1\}$. Recall that the product in reduced K-theory is given by 
$$\{A\}\cdot\{B\}=\{A\otimes B\}-{\rm rank}(A)\cdot \{B\}-{\rm rank}(B)\cdot \{A\}.$$ Therefore
$\big{(}a\{\gamma_1\}\big{)}\cdot\{\gamma_1\} = -2a\{\gamma_1\}$ and 
$$\{V\}\cdot\{\gamma_1\}  =  \{V\otimes\gamma_1\}-r\{\gamma_1\}-\{V\}= \{V\otimes\gamma_1\}-(r+a)\{\gamma_1\}.$$
Consequently, $\{V\otimes\gamma_1\}=(r-a)\{\gamma_1\}$ and assertion 1 follows. 
If in addition $V\cong V\otimes\gamma_1$, then $$2a\{\gamma_1\}=\{V\}+\{V\otimes\gamma_1\}=r\{\gamma_1\}$$
and we see from Theorem \ref{Adams} that  $2a\equiv r \mathrm{\ mod\ }2^{\phi(n)}$. \qed
\medbreak
Most of the vector bundles we will use in our study will be (isomorphic to) sub-bundles of the trivial vector bundle of rank $n+1$ over $\mathbb{R}P^n$. We will need the following observation about such vector bundles. 

\begin{lemma}
\label{subb}
Let $V$ be a vector bundle over $\mathbb{R}P^n$ which is isomorphic to a sub-bundle of the trivial vector bundle of rank $n+1$. Let $w(V)=p(x)\in  H^*(\mathbb{R}P^{n};\mathbb{Z}_2),$ where the degree of the polynomial $p$ is at most $n$. Then the degree of $p$ is either $0$ or $\mathrm{rank}(V)$.
\end{lemma}
\medbreak\noindent{\bf Proof.}
Let $W$ be a vector bundle over $\mathbb{R}P^n$ such that $V\oplus W$ is isomorphic to the trivial vector bundle of rank $n+1$ and let $w(W)=q(x)$. If $\mathrm{rank}(V)=r$, then the $\mathrm{rank}(W)=n+1-r$. It follows from property 1 of Stiefel-Whitney classes (see Theorem \ref{SW}) that the degree of $p(x)$ is at most $r$ and the degree of $q(x)$ is at most $n+1-r$. Therefore, the degree of $p(x)q(x)$ is at most $n+1$ with equality only in the case when the degree of $p(x)$ is $r$. Since the vector bundle $V\oplus W$ is trivial, we have: $$p(x)q(x)=1 \mathrm{\ in\ }H^*(\mathbb{R}P^{n};\mathbb{Z}_2).$$ Therefore, we either have $p(x)q(x)=1$, i.e. $p(x)=1$, or $p(x)q(x)=1+x^{n+1}$, i.e. the degree of $p(x)$ is equal to $r$.
\qed
\medbreak
We now put the two previous Lemmas together to obtain the main technical result we need in order to prove Wolf's Theorem. 
\begin{lemma}
\label{techn}
Let $V$ be a vector bundle over $\mathbb{R}P^n$ which is isomorphic to a sub-bundle of the trivial vector bundle of rank $n+1$. Let $k$ be an integer satisfying $\frac{n}{2}\le k\le n$ and let $\iota:\mathbb{R}P^k\to \mathbb{R}P^n$ be the natural inclusion.
\begin{enumerate}
\item If  the vector bundle $V\oplus (V\otimes\gamma_1)$ is isomorphic to a sub-bundle of the trivial vector bundle of rank $n+1$, then $\mathrm{rank}(V)=0$;
\item If $n\ne 1,3,7$ and if $V\cong\mathrm{Im}(S)$ for some continuous map $S:S^n\to \mathrm{Hom}(U,U)$ satisfying the conditions of Lemma \ref{oddmaps}, then $\mathrm{rank}(V)=0$;
\item If $k\ge 10$ and if $\iota^*(V)\cong\iota^*(V)\otimes\gamma_1$, then $\mathrm{rank}(V)=0$.
\end{enumerate}
\end{lemma}
\medbreak\noindent{\bf Proof.}
Let $\mathrm{rank}(V)=r$ and let $w(V)=p(x)\in H^*(\mathbb{R}P^{n};\mathbb{Z}_2)$. Assume that $V\oplus(V\otimes\gamma_1)$ is a sub-bundle of the trivial vector bundle of rank $n+1$; we then have $2r\le n+1$ and consequently $r\le n$. Let $w\big{(}V\oplus (V\otimes\gamma_1)\big{)}=q(x)\in H^*(\mathbb{R}P^{n};\mathbb{Z}_2).$
We see from Lemma \ref{subb} that the degree of $q(x)$ is either $0$ or $2r$. On the other hand, the inequality $r\le n$ and Lemma \ref{lema5} imply that $q$ is a polynomial of degree $r$; namely $q(x)=(1+x)^r$. Therefore, $r=0$.

Assume now $V\cong\mathrm{Im}(S)$ with $S:S^n\to \mathrm{Hom}(U,U)$ satisfying the conditions of Lemma \ref{oddmaps}. If $r\ne 0$, then $V$ allows a nowhere vanishing section and therefore the degree of $p$ is strictly smaller than $r$. It follows from Lemma \ref{subb} that $p(x)=1$. Since $V\cong V\otimes\gamma_1$, the relation $\{V\}+\{V\otimes\gamma_1\}=r\{\gamma_1\}$ implies 
$$(1+x)^r=1\mathrm{\ in\ } H^*(\mathbb{R}P^{n};\mathbb{Z}_2).$$
Hence if $r\ne0$, then it must be that  $r=n+1$, or equivalently $V$ is isomorphic to the trivial bundle of rank $n+1$. We use the isomorphism $V\cong V\otimes \gamma_1$ to conclude that $r\ne0$ implies $$(n+1)\{\gamma_1\}=0 \mathrm{\ in\ } \widetilde{KO}(\mathbb{R}P^n)$$
and in particular, $n+1\ge 2^{\phi(n)}$. Under our assumptions the last inequality is impossible (see Lemma \ref{phi}).

We now prove the last assertion of the Lemma. Let $j$ be the integer satisfying $2^j\le k<2^{j+1}$; the inequality $\frac{n}{2}\le k$ then implies $n<2^{j+2}$. Due to our assumption $k\ge 10$, we have $\phi(k)\ge j+3$; see Lemma \ref{phi}. Let $a$ be an integer such that $\{V\}=a\{\gamma_1\}$ in $\widetilde{KO}(\mathbb{R}P^n)$; then $\{\iota^*(V)\}=a\{\gamma_1\}$ in $\widetilde{KO}(\mathbb{R}P^k)$. It follows from Lemma \ref{lema5} that $2a\equiv r\mathrm{\ mod\ }2^{\phi(k)}$. Thus,
$$2a\equiv r\mathrm{\ mod\ }2^{j+3}\mathrm{\ \ and\ so\ \ } a\equiv \frac{r}{2}\mathrm{\ mod\ }2^{j+2}.$$
Note $r\le n+1$ and consequently $0\le\frac{r}{2}\le n<2^{j+2}$. Using $\{V\}=a\{\gamma_1\}$ we have $w(V)=(1+x)^a$ and therefore $$p(x)=(1+x)^{\frac{r}{2}}.$$ Hence the degree of the polynomial $p$ is $\frac{r}{2}$. Since by Lemma \ref{subb} the degree of the polynomial $p$ is either $0$ or $r$, we obtain $r=0$. 
\qed

\section{Vector Valued Polynomial Maps}

\hskip .2in  We now turn our attention to the polynomial aspect of our problem. 
 Let $V$ and $W$ be finite dimensional vector spaces and let $S^n(V)$ be the $n^{\mathrm{th}}$ symmetric 
power of $V$. A map $x:V\to W$ is said to be a {\it $W$-valued homogeneous polynomial of degree $n$} if there exists a linear map $x_L:S^n(V)\to W$ such that $x=x_L\circ i$, where $i:V\to S^n(V)$ is given by $v\mapsto v^n$. A map $P:V\to W$ is said to be a {\it $W$-valued polynomial} if it can be written as a sum of finitely many homogeneous polynomials.  The set of all polynomial maps $x:V\to W$ will be denoted by $\mathbb{R}[V,W]$. Note that $\mathbb{R}[V,\mathbb{R}]$ forms a ring under pointwise multiplication and that $\mathbb{R}[V,W]$ forms a module over $\mathbb{R}[V,\mathbb{R}]$.
\medbreak
The following are some of the polynomial maps we will use in our work.
\begin{itemize}
\item $Q:T_PM\to \mathbb{R}$ given by $Q(v)=(v,v)$, where $(.,.)$ denotes the metric of $M$ at the point $P$;
\item Elements $S\in\mathcal{P}$;
\item $x:T_PM\to T_PM$ given by $x(v)=S(v)x_0$, where $x_0\in T_PM$ and $S\in \mathcal{P}_n$ for some $n$.  
\end{itemize}

In general, if $V$ is equipped with a non-degenerate inner product of signature $(p,q)$, we view the inner product as a homogeneous element $Q\in\mathbb{R}[V,\mathbb{R}]$ of degree 2. Let $\mathcal{N}$ denote the set of all $v\in V$ with $Q(v)=0$ and let $\mathcal{I}$ denote the ideal in $\mathbb{R}[V,\mathbb{R}]$ generated by $Q$. 
Once an orthonormal basis for $V$ is chosen, we have isomorphisms $V\cong \mathbb{R}^{(p,q)}$ and $\mathbb{R}[V,\mathbb{R}]\cong\mathbb{R}[\lambda_1,\ldots,\lambda_{p+q}]$. Under these isomorphisms $Q$ corresponds to the polynomial
$$q=-\lambda_1^2-\ldots-\lambda_p^2+\lambda_{p+1}^2+\ldots+\lambda_{p+q}^2,$$
while $\mathcal{N}$ corresponds to the nullcone in $\mathbb{R}^{(p,q)}$, i.e. the nullset of the polynomial $q$. 

Note that the polynomial $z_1^2+z_2^2+z_3^2$ is irreducible over $\mathbb{C}$ and therefore, when $p+q>2$, the polynomial $q$ is irreducible over $\mathbb{C}$ as well. Hence, we may use Hilbert's Nullstellensatz in the proof of the following Lemma.

\begin{lemma}
\label{HN} 
Adopt the notation established above. 
Let $v\in \mathcal{N}$ and let $x:V\to W$ be a polynomial map. If $x=0$ on a neighborhood of $v$ in $\mathcal{N}$, then $x\in \mathcal{I}\mathbb{R}[V,W]$, i.e. there exists a polynomial map $y:V\to W$ with $x=Q\cdot y$.  Moreover, if $x$ is homogeneous of degree $n$, then $y$ is homogeneous of degree $n-2$.

\end{lemma}

\medbreak\noindent{\bf Proof.}
We shall, without loss of generality, assume that $v\ne0$. Let $\{e_1,\ldots, e_{p+q}\}$ be an orthonormal basis of $V$ and let $\{f_1,\ldots,f_w\}$ be a basis for $W$. These bases allow us to identify $V$ with the real part of $\mathbb{C}^{(p,q)}$, identify $W$ with the real part of $\mathbb{C}^w$, and identify $\mathcal{N}$ with the real part of the complex nullcone 
$$\mathcal{N}^\mathbb{C}:=\{(z_1,\ldots, z_{p+q})|-z_1^2-\ldots-z_p^2+z_{p+1}^2+\ldots+z_{p+q}^2=0\}.$$
The polynomial map $x:V\to W$ can now be considered as a collection $x_1,\ x_2,\ \ldots,\ x_w$ of $w$ polynomials in $p+q$ variables over $\mathbb{C}$. These polynomials, by our assumption, vanish on a neighborhood of $v$ in $\mathcal{N}^{\mathbb{C}}\cap\mathbb{R}^{(p,q)}$.

Let $v_i$,  where $1\le i\le p$, be a non-zero coordinate of $v$. Consider a branch of $\sqrt{z}$ around $v_i$ and a holomorphic chart 
$$\phi: (z_1,\ldots, \widehat{z_i}, \ldots, z_{p+q})\mapsto (z_1,\ldots,\sqrt{-z_1^2-\ldots-\widehat{z_i^2}-\ldots+z_{p+q}^2}\ ,\ldots, z_{p+q})$$ 
around $v$. The holomorphic functions  $x_i\circ \phi$ each vanish on a real neighborhood of  the point $(v_1,\ldots,\widehat{v_i},\ldots,v_{p+q})$. It now follows from the Identity Theorem that the functions $x_i\circ \phi$ vanish on a complex neighborhood of the point $(v_1,\ldots,\widehat{v_i},\ldots,v_{p+q})$.  In other words, polynomials $x_i$ vanish on a (complex) neighborhood of $v$.  By analytic continuation, the  holomorphic functions $x_i$ vanish on the entire $\mathcal{N}^{\mathbb{C}}$.
 
Since the polynomials $x_i$ vanish on the nullset of the polynomial 
$$q=-z_1^2-\ldots-z_p^2+z_{p+1}^2+\ldots+z_{p+q}^2,$$ Hilbert's Nullstellensatz implies the polynomials $x_i$ are divisable by $q$.
Thus $x_i=q\cdot y_i$ for some polynomials $y_i$. Note that since polynomials $x_i$ and $q$ have real coefficients, so do polynomials $y_i$. It is now immediate that there exists a polynomial map $y:V\to W$ such that $x=Q\cdot y.\qed$ 
\medbreak
Let $x_1,\ldots,x_r:V\to W$ be polynomial maps. Assume that $r\le w$, where $w=\dim W$. Our next step is to relate linear dependence of vectors $x_1(v),\ldots,x_r(v)\in W$ with linear dependence of polynomial maps $x_1,\ldots,x_r$ over the ring $\mathbb{R}[V,\mathbb{R}]$.
\medbreak
Let $\{f_1,\ldots, f_w\}$ be a basis for $W$. A collection of $r$ polynomial maps $x_1,\ldots,x_r:V\to W$ can be identified with a $w\times r$ matrix $\mathbf{X}$ with entries in $\mathbb{R}[V,\mathbb{R}]$.  Let $$\big{\{}M_1,\ldots, M_{w\choose r}\big{\}}\subset \mathbb{R}[V,\mathbb{R}]$$ be the set of all $r\times r$ minors of $\mathbf{X}$. Let $I(x_1,\ldots, x_r)$ denote the ideal in $\mathbb{R}[V,\mathbb{R}]$ generated by the minors $M_1,\ldots, M_{w\choose r}$. We have the following Lemma. 

\begin{lemma}
Adopt the notation established above and assume $r\le w$. The ideal $I(x_1,\ldots,x_r)$ of $\mathbb{R}[V,\mathbb{R}]$ is independent of the choice of the basis $\{f_1,\ldots, f_w\}$ of $W$. 
\end{lemma}

\medbreak\noindent{\bf Proof.}
Changing the basis of $W$ amounts to multiplication of $\mathbf{X}$ on the left by an element of $GL(w,\mathbb{R})$. The Lemma now follows from the fact that elementary row operations (over $\mathbb{R}$) have one of the following three effects on $\mathbf{X}$:
\begin{itemize}
\item Preserve the set of generators $\big{\{}M_1,\ldots, M_{w\choose r}\big{\}}$;
\item Replace certain minors $M_i$ with $k\cdot M_i$ for some real number $k\ne 0$;
\item Replace certain minors $M_i$ with $M_i+k\cdot M_j$ for some real number $k$ and $i\ne j$.\qed
\end{itemize}
\medbreak
We would like to point out the following properties of the ideal $I(x_1,\ldots, x_r)$; they are immediate consequences of the corresponding properties of the minors $M_i$ of the matrix $\mathbf{X}$.

\begin{lemma}
\label{prop_I}
Adopt the notation established above and assume $r\le w$.
\begin{enumerate}
\item If $\sigma$ is any permutation of indices $\{1,\ldots,r\}$, then $I(x_{\sigma_1},\ldots,x_{\sigma_r})=I(x_1,\ldots,x_r)$;
\item If $c\in\mathbb{R}[V,\mathbb{R}]$, then $I(cx_1, x_2, \ldots,x_r)=(c^r)I(x_1, x_2, \ldots,x_r),$ where $(c^r)$ is the ideal of the ring $\mathbb{R}[V,\mathbb{R}]$ generated by $c^r$;
\item If $c\in \mathbb{R}[V,\mathbb{R}]$, then $I(x_1+cx_2,x_2,\ldots,x_r)=I(x_1,x_2,\ldots,x_r).\qed$
\end{enumerate}
\end{lemma}

We use the ideal $I(x_1,\ldots,x_r)$ as a tool for studying linear dependence of vectors $x_1(v),\ldots, x_r(v)$.

\begin{lemma}
\label{ind1}
Adopt the notation established above and assume $r\le w$. The following are equivalent:
\begin{enumerate}
\item Vectors $x_1(v),\ldots, x_r(v)\in W$ are linearly dependent for all $v\in V$;
\item We have $I(x_1,\ldots, x_r)=\{0\}$.\qed
\end{enumerate}
\end{lemma}

Note that $\bigcap_{k\in\mathbb{N}}\mathcal{I}^k=\{0\}.$ Therefore, for a given collection $x_1,\ldots, x_r\in\mathbb{R}[V,W]$ we either have $$I(x_1,\ldots,x_r)=\{0\} \mathrm{\  or\ } I(x_1,\ldots,x_r)\not\subset \mathcal{I}^{k+1} \mathrm{\ for\  some\ } k\ge 0.$$
In the case when $I(x_1,\ldots,x_r)\ne\{0\}$ we let $k(x_1,\ldots,x_r)$ denote the smallest $k\in\mathbb{N}\cup\{0\}$ such that $I(x_1,\ldots, x_r)\not\subset\mathcal{I}^{k+1}$. We will refer to $k(x_1,\ldots, x_r)$ as {\it the degree of the linear dependence of $x_1,\ldots, x_r$ over the nullcone}. 

\begin{lemma}
\label{ind2}
Adopt the notation established above and assume $r\le w$. The following are equivalent:
\begin{enumerate}
\item Vectors $x_1(v),\ldots, x_r(v)\in W$ are linearly dependent for all $v\in \mathcal{N}$;
\item We have $I(x_1,\ldots, x_r)\subset\mathcal{I}$;
\item There exist $c_1,c_2,\ldots,c_r\in\mathbb{R}[V,\mathbb{R}]$ with $c_i\not\in\mathcal{I}$ for at least one $i$ and $$c_1x_1+\ldots+c_rx_r\in\mathcal{I}\mathbb{R}[V,W].$$
\end{enumerate}
\end{lemma}

\medbreak\noindent{\bf Proof.}
It is clear that assertion 2 implies assertion 1. To justify that 1 implies 2, we use Lemma \ref{HN}, which says that the only elements of $\mathbb{R}[V,\mathbb{R}]$ vanishing over $\mathcal{N}$ are those in $\mathcal{I}$. 

We now prove that the third assertion implies the first. Suppose vectors $x_1(v_0),\ldots, x_r(v_0)$ are linearly independent for some $v_0\in\mathcal{N}$; then the vectors $x_1(v),\ldots, x_r(v)$ are linearly independent for all $v$ in a neighborhood $\mathcal{B}\subset\mathcal{N}$ of $v_0\in\mathcal{N}$. Since
$$c_1(v)x_1(v)+\ldots+c_r(v)x_r(v)=0 \mathrm{\ \ for\ all\ \ } v\in\mathcal{B},$$
we have $c_1(v)=\ldots=c_r(v)=0$ for all $v\in\mathcal{B}$. It now follows from Lemma \ref{HN} that $c_1,\ldots, c_r\in\mathcal{I}$. Contradiction.

We now show 1 implies 3. Let $q(v)$ be the size of a maximal linearly independent subset of $x_1(v),\ldots,x_r(v)$. Let $q=\max_{v\in\mathcal{N}}q(v)$; by assumption $q<r$. Let $v_0$ be such that $q(v_0)=q$. We may, without loss of generality, assume that the vectors $x_1(v_0),\ldots,x_q(v_0)$ are linearly independent. Then there exists a $q\times q$ minor, $M\in\mathbb{R}[V,\mathbb{R}]$, of the matrix $\widetilde{\mathbf{X}}$ corresponding to the collection $x_1,\ldots,x_q$  such that $M(v_0)\ne 0$ and consequently $M\not\in\mathcal{I}$. Without loss of generality we may assume the basis for $W$ is chosen so that  the minor $M$ comes from the submatrix of $\widetilde{\mathbf{X}}$ consisting of first $q$ rows.  

Now consider polynomials $x_1,\ldots,x_q,x_r$. By the definition of $q$, vectors $x_1(v),\ldots,x_q(v),x_r(v)$ are linearly dependent for all $v\in\mathcal{N}$. Since $x_1(v),\ldots,x_q(v)$ are linearly independent over a neighborhood $\mathcal{B}\subset\mathcal{N}$ of $v_0$ (for example, $\mathcal{B}$ with $M(v)\ne 0$ for all $v\in\mathcal{B}$), there exist coefficients $c_1(v),\ldots,c_q(v)$ with 
\begin{equation}
\label{lin-ind}
x_r(v)=c_1(v)x_1(v)+\ldots+c_q(v)x_q(v) \mathrm{\ \ for\ all\ \ } v\in\mathcal{B}.
\end{equation}
The coefficients $c_i(v)$ depend rationally on $v$. Indeed, they are determined by a $q\times q$ system 
of linear equations with coefficients in $\mathbb{R}[V,\mathbb{R}]$; this system comes from considering the first $q$ coordiantes of the vectors on both sides of (\ref{lin-ind}). Note that the determinant of this system is equal to $M$. Therefore, $$c_i(v)=\frac{p_i(v)}{M(v)},\ \ \mathrm{with}\ \ v\in\mathcal{B}\  \mathrm{and}\ p_i,M\in\mathbb{R}[V,\mathbb{R}].$$
We now see from display (\ref{lin-ind}) that the identity   
$$p_1(v)x_1(v)+\ldots+p_q(v)x_q(v)+(-M(v))x_r(v)=0$$
holds for all $v\in\mathcal{B}$. 
It follows from Lemma \ref{HN} that $$p_1x_1+\ldots+p_qx_q+(-M)x_r\in\mathcal{I}\mathbb{R}[V,\mathbb{R}]. \ \qed$$ 

The following result will play a crucial role in our proof of Wolf's Theorem. 

\begin{lemma}
\label{ind3}
Adopt the notation established above and assume $r\le w$.
Let $\Delta$ be a $\mathbb{R}[V,\mathbb{R}]$-submodule of $\mathbb{R}[V,W]$ such that  
\begin{equation}
\label{cond}
Q\cdot x\in \Delta \Longrightarrow x\in \Delta.
\end{equation}
If there exist $x_1,\ldots,x_r\in\Delta$ with $I(x_1,\ldots, x_r)\ne\{0\}$, then there exist $x_1,\ldots, x_r\in\Delta$ with $k(x_1,\ldots, x_r)=0$.
\end{lemma}

\medbreak\noindent{\bf Proof.}
Let  $$\mathcal{A}=\Big{\{}(x_1,\ldots,x_r)\Big{|}x_1,\ldots,x_r\in \Delta,\ I(x_1,\ldots,x_r)\ne\{0\}\Big{\}}$$ and let $s=\min_{(x_1,\ldots,x_r)\in\mathcal{A}} k(x_1,\ldots,x_r)$. We have $I(x_1,\ldots,x_r)\subset \mathcal{I}^s$ for all $(x_1,\ldots,x_r)\in\mathcal{A}$. Assume now $s\ge 1$. Let $x_1,\ldots,x_r$ be such that $k(x_1,\ldots,x_r)=s$; in particular 
$$I(x_1,\ldots,x_r)\not\subset\mathcal{I}^{s+1}\mathrm{\ \ and\ \ } I(x_1,\ldots,x_r)\subset \mathcal{I}.$$
By the previous Lemma there exist $c_1,\ldots,c_r\in\mathbb{R}[V,\mathbb{R}]$ such that  $c_i\not\in\mathcal{I}$ for some $i$ and $$c_1x_1+\ldots+c_rx_r\in\mathcal{I}\mathbb{R}[V,\mathbb{R}].$$
Without loss of generality we may assume $i=1$.

Let $y\in\mathbb{R}[V,\mathbb{R}]$ be such that $c_1x_1+\ldots+c_rx_r=Q\cdot y$. Note that $Q\cdot y\in \Delta$ and consequently $y\in \Delta$. We now study the collection $(y,x_2,\ldots,x_r)$. The properties of Lemma \ref{prop_I} imply that 
\begin{eqnarray}
\label{ideals}
&(Q^r)I(y,x_2,\ldots,x_r)=I(Qy,x_2,\ldots,x_r)=I(c_1x_1+\ldots+c_rx_r,x_2,\ldots,x_r)\\
&\ \ \ \ \ \ \ \ \ \ \ \ \ \ \ =I(c_1x_1,x_2,\ldots,x_r)=(c_1^r)I(x_1,x_2,\ldots,x_r).\label{ideals2}
\end{eqnarray}
Since $I(x_1,\ldots,x_r)\ne\{0\}$, we have $I(y,x_2,\ldots,x_r)\ne\{0\}$. Therefore $(y,x_2,\ldots,x_r)\in\mathcal{A}$ and, by definition of $s$, $I(y,x_2,\ldots,x_r)\subset \mathcal{I}^s$. Relations of display (\ref{ideals}-\ref{ideals2}) now imply 
$$(c_1^r)I(x_1,\ldots,x_r)\subset \mathcal{I}^{s+r}\subset\mathcal{I}^{s+1}.$$
We now see that the polynomial map $Q^{s+1}$ divides $c_1^rP$ for all $P\in I(x_1,\ldots,x_r)$. Since  $Q$ is a prime element of $\mathbb{R}[V,\mathbb{R}]$ and since $c_1$ is not divisible by $Q$ (by assumption $c_1\not\in \mathcal{I}$), the element $Q^{s+1}$ must divide $P$ for all $P\in I(x_1,\ldots,x_r)$. This contradicts our assumption that $I(x_1,\ldots,x_r)\not \subset \mathcal{I}^{s+1}$.\qed

\section{Signature $(p,p)$}

\hskip .2in In this Section we complete the proof of  Wolf's Theorem for locally isotropic manifolds of signature $(p,p)$, where $p\ne 2,4,8$. As pointed out in the Introduction, nilpotency of the Szab\'o operator over the nullcone plays a big role in our proof. We start the proof by extending the nilpotency result of Gilkey-Stavrov \cite{GilkeyStavrov}.

\begin{theorem}
\label{thm5}
Let $M$ be a locally isotropic manifold and let $P\in M$.  If $S\in\mathcal{P}$, then $S(v)$ is nilpotent for all $v\in \mathcal{N}$.
\end{theorem}

\medbreak\noindent{\bf Proof.}
There is nothing to show in the case $v=0$.
When  $v\in \mathcal{N}-\{0\}$, take $\lambda\in \mathbb{R}$ with $\lambda\ne 0$. We have $(v,v)=(\lambda v,\lambda v)=0$. Since manifold $M$ is locally isotropic, there exists a local isometry $\mathcal{T}$ fixing $P$ such that its differential $T:T_PM\to T_PM$ satisfies
$$
T(v)=\lambda v, \ \  S(\lambda v)=T\circ S(v)\circ T^{-1}.
$$
Consequently, 
$$
{\rm Tr}\{S(\lambda v)^n\}={\rm Tr}\{S(v)^n\}{\rm \ \ for\ \ all\ \ } n\in \mathbb{N}.
$$
Since $\lambda$ was arbitrary, we may take the limit as $\lambda\to0$. As the map $S:T_PM\to \mathrm{Hom}(T_PM,T_PM)$ is continuous, we obtain
${\rm Tr}\{S(v)^n\}=0$ for all $n\in \mathbb{N}. $\qed
\medbreak
The simpliest case of nilpotency over the nullcone is if an operator $S\in\mathcal{P}$ vanishes over the nullcone. 
\begin{lemma}
\label{nullstellensatz}
Let $M$ be a pseudo-Riemannian manifold and let $P\in M$. Let $S\in\mathcal{P}_{P}$ be such that $S\equiv 0$ over the nullcone $\mathcal{N}\subset T_PM$. Then:
\begin{enumerate}
\item There exists $T\in\mathcal{P}_{P}$ such that $S(v)=(v,v)T(v)$; 
\item $S\in\mathcal{P}_{1}$ implies $S\equiv 0$ over $T_PM$. 
\end{enumerate}
\end{lemma}
\medbreak\noindent{\bf Proof.}
The first assertion follows as an immediate consequence of Lemma \ref{HN}. To prove the second claim it is enough to prove that $\mathcal{P}_0=\{0\}$.

Let $S\in \mathcal{P}_0$. The self-adjoint operator $S(v)$ depends linearly upon $v$ and satisfies $S(v)v=0$ for all $v\in T_PM$. Polarizing the last identity we get $S(v)w+S(w)v=0$ for all $v$ and $w$. Consequently:
$$
0=(S(v)w,w)+(S(w)v,w)=(S(v)w,w)+(v,S(w)w)=(S(v)w,w).
$$
Further polarization yields $(S(v)x,y)+(x,S(v)y)=0$, i.e. $(S(v)x,y)=0$ for all $x,y$ and $v$. Thus, the operator $S(v)$ is the zero operator for all $v$ and hence $S\equiv 0.\qed$
\medbreak
The following step is to study nilpotency of order 
$2$,  i.e. the case when $S^2\equiv 0$ over the nullcone. 
   
\begin{lemma}
\label{lema1}
Let  $M$ be a locally isotropic pseudo-Riemannian manifold of signature $(p,p)$, where $p\ne 2,4,8$.
Let $S\in \mathcal{P}_P$ be such that $S^2\equiv 0$ over the nullcone bundle. Then:
\begin{enumerate}
\item $S\equiv 0$ over $\mathcal{N}$; 
\item $S\in\mathcal{P}_{1}$ implies $S\equiv 0$ over $T_PM$. 
\end{enumerate}
\end{lemma}
\medbreak\noindent{\bf Proof.}
The tangent space $T_PM$ decomposes into the direct sum of a maximal negative definite subspace $V_-$ and its orthogonal complement $V_+$. Let
$\varrho_+:V\to V_+$ and $\varrho_-:V\to V_-$ denote the corresponding
orthogonal projections and let $\Phi:V\to V$ be the linear map given by
$\Phi :=\varrho_+ -\varrho_- .$ The decomposition $V_-\oplus V_+$ gives rise to an embedding $i:S^{p-1}\to\mathcal{N}$; namely, $i(v)=(v,v)$. 

Consider the map $\tilde{S}=S\circ i:S^{p-1}\to \mathrm{Hom}(T_PM,T_PM)$. This map satisfies the conditions of Lemma \ref{oddmaps} and thus gives rise to a vector bundle $\mathrm{Im}(\tilde{S})$ over $\mathbb{R}P^{p-1}$. Let $r$ be the rank of $\mathrm{Im}(\tilde{S})$. To prove assertion 1 it is enough to prove $r=0$. We assume  $r>0$ and argue for contradiction. 

The subspaces ${\rm Im}(\tilde{S}(v))$ are totally isotropic for all $v\in S^{p-1}$. Indeed,
we may use the self-adjointness of $\tilde{S}$ to compute:
$$(\tilde{S}(v)x,\tilde{S}(v)y)=(x,\tilde{S}^2(v)y)=0 {\rm\ \  for\ \  all\ \  } x,y\in T_PM.$$
Totally isotropic subspaces of $T_PM$ project isomorphically into the negative
definite subspace $V_-$ via the orthogonal projection.
Thus, the map $$\varrho_-:{\rm Im}(\tilde{S})\to \mathbb{R}P^{p-1}\times V_-,$$
defined by  the projection $\varrho_-$ on fibers, maps ${\rm Im}(\tilde{S})$ isomorphically onto a sub-bundle of the trivial bundle of rank $p$. However, by Lemma \ref{techn} we see that ${\rm Im}(\tilde{S})$ is not isomorphic to a sub-bundle of the trivial bundle of rank $p$. This contradiction completes the proof of the first statement. The second statement follows as an immediate corollary to the previous Lemma. \qed
\medbreak
One of the consequences of Lemma \ref{lema1} is the following result.

\begin{lemma}
\label{lema3}
Let  $M$ be a locally isotropic pseudo-Riemannian manifold of signature $(p,p)$, where $p\ne 2,4,8$.
Let $S\in \mathcal{P}$. Then $S^3\equiv 0$ over the nullcone.
\end{lemma}

\medbreak\noindent{\bf Proof.}
By Theorem \ref{thm5} we know that $S$ is nilpotent over the nullcone. Let $n$ be the smallest integer such that $S^n \equiv 0$ over the nullcone. If $n\le 3$ there is nothing to show. So, we assume $n\ge 4$ and argue for contradiction. Let $k$ be the greatest odd number smaller than $n$, i.e. let
$$k=\begin{cases}
      {n-1& {\rm\  for\  $n$\  even}, \cr
      n-2 & {\rm\  for\  $n$\  odd}.}
\end{cases}$$
By the choice of $n$ we have that $T:=S^k\not\equiv 0$ over the nullcone. Since $n\ge 4$ implies $2k\ge n$, we have $T^2\equiv 0$ over the nullcone. Note also that $T\in \mathcal{P}$. Thus, by Lemma \ref{lema1}, we have that $T\equiv 0$ over the nullcone. This contradiction completes the proof of the Lemma.\qed
\medbreak
We now establish a result from linear algebra regarding self-adjoint operators $A$ on vector spaces of indefinite signature satisfying  $A^3\equiv 0$.

\begin{lemma}
\label{lema4}
Let $V$ be a vector space with non-degenerate inner product $(.,.)$ of signature $(p,q)$. Decompose $V=V_-\oplus V_+$ as a direct sum of a positive definite subspace $V_+$ and its negative definite orthogonal complement $V_-$; let
$\varrho_+:V\to V_+$ and $\varrho_-:V\to V_-$ denote the corresponding
orthogonal projections. Define a linear map $\Phi:V\to V$ by
$\Phi v:=\varrho_+ v-\varrho_- v.$
Let $A$ be a self-adjoint map on $V$ such that $A^3=0$. We have:
\begin{enumerate}
\item The subspace ${\rm Im}(A^2)$ is totally isotropic;
\item ${\rm rank}(A\Phi A^2)={\rm rank}(A^2\Phi A^2)={\rm rank}(A^2)$ and the map $A:{\rm Im}(A\Phi A^2)\to {\rm Im}(A^2)$ is an isomorphism;
\item The subspace ${\rm Im}(A\Phi A^2)$ inherits a non-degenerate inner product;
\item Subspaces ${\rm Im}(A^2)$ and ${\rm Im}(A\Phi A^2)$ are orthogonal with respect to $(.,.)$.
\end{enumerate}
\end{lemma}

\medbreak\noindent{\bf Proof.}
Assertion 1 is immediate from the assumptions that $A$ is self-adjoint with $A^3=0$:
$$(A^2x,A^2y)=(Ax,A^3y)=0 {\rm\ \  for\ \  all\ \  } x,y\in V.$$
To prove assertion 2 we first consider the map $A^2\Phi$ - it is of the same rank as $A^2$ and it is self-adjoint with respect to the positive definite inner product
$g$ given by
$g(v,w):=(v,\Phi w)=(\Phi v,w).$ Thus, the map $A^2\Phi$ is an automorphism of ${\rm Im}(A^2\Phi)={\rm Im}(A^2)$. Consequently, 
$${\rm rank}(A^2)\ge {\rm rank (A\Phi A^2)}\ge{\rm rank}(A^2\Phi A^2)={\rm rank}(A^2),$$
and the desired statements follow. 

Since $A^2\Phi$ is an isomorphism of $\mathrm{Im}(A^2)$, we have $g(A^2\Phi x, A^2\Phi x)\ge 0$ for all non-zero $x\in\mathrm{Im}(A^2)$. In other words, we have $g(A^2\Phi A^2x, A^2\Phi A^2x)\ge 0$ with equality if and only if $A^2x=0$. Let $x$ be such that $A\Phi A^2x\ne 0$, or equivalently $A^2x\ne 0$. To prove that $\mathrm{Im}(A\Phi A^2)$ inherits the non-degenerate inner product, it is enough to find $y$ such that $(A\Phi A^2x, A\Phi A^2y)\ne 0$. Consider $y=\Phi A^2x$. We see from $A^2x\ne 0$ that
$$(A\Phi A^2x, A\Phi A^2y)=(A\Phi A^2x,A\Phi A^2\Phi A^2x)=g(A^2\Phi A^2x, A^2\Phi A^2 x)\ne0.$$ 
Finally, we verify assertion 4 by computing:
$(A^2x, A\Phi A^2y)=(A^3x, \Phi A^2y)=0.\qed$
\bigbreak
We are now ready to prove Wolf's Theorem in signatures $(p,p)$ where $p\ne 2,4,8.$

\begin{theorem}
Let  $M$ be a locally isotropic pseudo-Riemannian manifold of signature $(p,p)$, where $p\ne 2,4,8$.
Then $M$ is locally symmetric. 
\end{theorem}

\medbreak\noindent{\bf Proof.}
Fix a point $P\in M$. The tangent space $T_PM$ decomposes into the direct sum of a maximal negative definite subspace $V_-$ and its orthogonal complement $V_+$. Let
$\varrho_+:V\to V_+$ and $\varrho_-:V\to V_-$ denote the corresponding
orthogonal projections and  let a linear map $\Phi:V\to V$ be defined by
$\Phi v:=\varrho_+ v-\varrho_- v.$ The decomposition $V_-\oplus V_+$ gives rise to an embedding $i:S^{p-1}\to\mathcal{N}$; namely, $i(v)=(v,v)$. 

We see from the previous two lemmas that the Szab\'o operator $\mathbf{S}(v)$ satisfies  
$${\rm rank}(\mathbf{S}(v)\Phi \mathbf{S}(v)^2)={\rm rank}(\mathbf{S}(v)^2)={\rm const\ \ for\ all\ }v\in\mathcal{N}-\{0\}.$$
In particular, this means that there exist vector bundles $E$ and $F$ over the nullcone $\mathcal{N}-\{0\}$ having fibers 
$$E\big{|}_v:=\mathrm{Im}\big{(}\mathbf{S}(v)^2\big{)}\mathrm{\ \ and\ \ }
F\big{|}_v:=\mathrm{Im}\big{(}\mathbf{S}(v)\Phi \mathbf{S}(v)^2\big{)}.$$
These vector bundles are sub-bundles of the trivial vector bundle  $\big{(}\mathcal{N}-\{0\}\big{)}\times V$. 

We now consider the vector bundles $i^*(E)$ and $i^*(F)$. From part 2 of the previous Lemma we see that the vector bundle map $\mathcal{S}:i^*(F)\to i^*(E)$ defined on fibers by the map $\mathbf{S}(i(v))$ is a vector bundle isomorphism. Note that both $i^*(E)$ and $i^*(F)$ descend to define vector bundles over $\mathbb{R}P^{p-1}$, which we will denote by $\mathrm{Im}(\mathbf{S}^2)$ and $\mathrm{Im}(\mathbf{S}\Phi\mathbf{S}^2)$ respectively. By property $\mathbf{S}(i(-v))=-\mathbf{S}(i(v))$ the vector bundle isomorphism $\mathcal{S}$ does not descend to an isomorphism between $\mathrm{Im}(\mathbf{S}\Phi\mathbf{S}^2)$ and $\mathrm{Im}(\mathbf{S}^2)$. Rather, it gives rise to a vector bundle isomorphism $\mathrm{Im}(\mathbf{S}\Phi\mathbf{S}^2)\cong \mathrm{Im}(\mathbf{S}^2)\otimes\gamma_1$.

We see from Lemma \ref{lema4} that $\mathrm{Im}(\mathbf{S}\Phi\mathbf{S}^2)$ and $\mathrm{Im}(\mathbf{S}^2)$ are fiber-wise orthogonal sub-bundles of the trivial vector bundle $\mathbb{R}P^{p-1}\times V$, that $\mathrm{Im}(\mathbf{S}\Phi\mathbf{S}^2)$ inherits non-degenerate metric and that $\mathrm{Im}(\mathbf{S}^2)$ is totally isotropic.
Thus, we may decompose:
$$\mathrm{Im}(\mathbf{S}\Phi\mathbf{S}^2)= F_-\oplus F_+,$$
where $F_-$ is a maximal negative definite sub-bundle and $F_+$ is its orthogonal complement. Let $\mathbf{S}F_-$ be the sub-bundle of $\mathrm{Im}(\mathbf{S}^2)$ over $\mathbb{R}P^{p-1}$ corresponding to $F_-\otimes\gamma_1$ under the isomorphism $\mathrm{Im}(\mathbf{S}^2)\cong\mathrm{Im}(\mathbf{S}\Phi\mathbf{S}^2)\otimes\gamma_1$; we have $F_-\otimes \gamma_1\cong \mathbf{S}F_-$.

We now study the vector bundle $F_-\oplus \mathbf{S}F_-$ over $\mathbb{R}P^{p-1}$.  We may apply Lemma \ref{techn} to see that either ${\rm rank}(F_-)=0$ or the vector bundle $F_-\oplus \mathbf{S}F_-$ is not isomorphic to a sub-bundle of the trivial bundle of rank $p$. Since  $F_-$ is negative definite and fiber-wise orthogonal to the totally isotropic $\mathbf{S}F_-$, no fiber of $F_-\oplus \mathbf{S}F_-$ contains a spacelike vector. Thus, the orthogonal projection $\varrho_-:V\to V_-$ gives rise to an isomorphism between $F_-\oplus \mathbf{S}F_-$ and a sub-bundle of $\mathbb{R}P^{p-1}\times V_-$. 
It now follows that ${\rm rank}(F_-)=0$ and that $\mathrm{Im}(\mathbf{S}\Phi\mathbf{S}^2)$ inherits a positive definite inner product.

We now consider $\mathrm{Im}(\mathbf{S}\Phi\mathbf{S}^2)\oplus \mathrm{Im}(\mathbf{S}^2)$ and apply the same argument as above. Since we have an isomorphism  
$\mathrm{Im}(\mathbf{S}\Phi\mathbf{S}^2)\otimes\gamma_1\cong \mathrm{Im}(\mathbf{S}^2)$, we may apply Lemma \ref{techn}. It follows that either ${\rm rank}(\mathrm{Im}(\mathbf{S}\Phi\mathbf{S}^2))={\rm rank}(\mathrm{Im}(\mathbf{S}^2))=0$ or the vector bundle $\mathrm{Im}(\mathbf{S}\Phi\mathbf{S}^2)\oplus \mathrm{Im}(\mathbf{S}^2)$ is not isomorphic to a sub-bundle of the trivial bundle of rank $p$. Since  $\mathrm{Im}(\mathbf{S}\Phi\mathbf{S}^2)$ is positive definite and fiber-wise orthogonal to the totally isotropic $\mathrm{Im}(\mathbf{S}^2)$, no fiber of $\mathrm{Im}(\mathbf{S}\Phi\mathbf{S}^2)\oplus \mathrm{Im}(\mathbf{S}^2)$ contains a timelike vector. Thus, the orthogonal projection $\varrho_+:V\to V_+$ gives rise to an isomorphism between $\mathrm{Im}(\mathbf{S}\Phi\mathbf{S}^2)\oplus \mathrm{Im}(\mathbf{S}^2)$ and a sub-bundle of $\mathbb{R}P^{p-1}\times V_+$. 
It now follows that ${\rm rank}(\mathrm{Im}(\mathbf{S}^2))=0$ and that $\mathbf{S}^2\equiv 0$ over the nullcone. 

The Theorem is now an immediate corollary to Lemma \ref{lema1} and Theorem \ref{thm2}.\qed

\section{The General Case}

\hskip .2in Let $P$ be a point of a locally isotropic pseudo-Riemannian manifold $M$ of signature $(p,q)$. In this Section we assume $q>p\ge 2$; the corresponding results in the case $p>q$ follow from the ones in the case $q>p$ be reversing the sign of the inner product. Let  $x$ and $y$ be two unit spacelike  (two unit timelike or two non-zero null) vectors at $P$ and let $S\in\mathcal{P}$. Since $M$ is locally isotropic, there exists $T:T_PM\to T_PM$ such that 
$$S(y)=T\circ S(x)\circ T^{-1}.$$
Thus, the rank, the spectrum, the Jordan normal form and the minimal polynomial of the operator $S(v)$ are all independent of the choice of unit spacelike (resp. timelike or non-zero null) vector $v$.
Let ${\rm Spec}^+(S)$ (resp. ${\rm Spec}^-(S)$) denote the spectrum of the operator $S$ over the unit spacelike (resp. timelike) pseudo-sphere in $T_PM$. Let $r_{+}$ (resp. $r_{-}$, $r_{0}$) denote the rank of $S$ over the unit spacelike pseudo-sphere (resp. timelike pseudo-sphere, $\mathcal{N}-\{0\}$) in $T_PM$. These satisfy the following relations (see \cite{GilkeyReniJa}).
\begin{lemma}
\label{thm4}
Adopt the notation established above and assume $S\in \mathcal{P}_{n}$ for some $n\ge 1$. We have:
\begin{enumerate}
\item $r_{-}=r_{+}$;
\item If $S\not\equiv 0$, then $r_0<r_+$.
\end{enumerate}
\end{lemma}

From this point on we simplify the notation by setting $r=r_-=r_+$. Note that the relation $E\big{|}_v:=\mathrm{Im}\big{(}S(v)\big{)}$ defines a vector bundle over $T_PM-\mathcal{N}$ of rank $r$. But, if $S\not\equiv 0$, the relation $E\big{|}_v=\mathrm{Im}\big{(}S(v)\big{)}$ does {\it not} define a vector bundle over $T_PM-\{0\}\simeq S^{p+q-1}$. 
What we are about to do is construct a vector bundle $E$ over $T_PM$ such that $E\big{|}_v=\mathrm{Im}\big{(}S(v)\big{)}$ for $v\not\in\mathcal{N}$.

Let $\iota_1:\mathbb{R}^{p}\to T_PM$ be an inclusion of $\mathbb{R}^p$ as a negative definite subspace  and let $\iota_2:\mathbb{R}^{q}\to T_PM$ be an inclusion of $\mathbb{R}^q$ as a positive definite subspace of $T_PM$.  Consider the  continuous maps
$$S\circ\iota_1:S^{p-1}\to \mathrm{Hom}(T_PM,T_PM) \mathrm{\ \ and\ \ } S\circ\iota_2:S^{q-1}\to \mathrm{Hom}(T_PM,T_PM).$$
The maps $S\circ\iota_1$ and $S\circ\iota_2$ satisfy the conditions of Lemma \ref{oddmaps} and we therefore have vector bundles $E_1$ and $E_2$ over $S^{p-1}$ and $S^{q-1}$, respectively. These vector bundles satisfy 
\begin{equation}
\label{e1e2}
E_1\big{|}_v:=\mathrm{Im}\big{(}S(\iota_1v)\big{)} \mathrm{\ \ and\ \ } E_2\big{|}_v:=\mathrm{Im}\big{(}S(\iota_2v)\big{)}.
\end{equation}
Moreover, they descend to vector bundles $\mathrm{Im}(S)_1$ and $\mathrm{Im}(S)_2$ over $\mathbb{R}P^{p-1}$ and $\mathbb{R}P^{q-1}$, respectively. We know that 
$\mathrm{Im}(S)_i\cong\mathrm{Im}(S)_i\otimes\gamma_1.$

To "glue" vector bundles $E_1$ and $E_2$ over the nullcone, we need to take a more global approach and use vector valued polynomial maps. The starting point of our construction is the following Lemma, also proved in \cite{GilkeyReniJa}.

\begin{lemma}
\label{thm4.5}
Adopt the notation established above and assume $S\in \mathcal{P}_{n}$ for some $n\ge 1$. We have:
\begin{enumerate}
\item ${\rm Spec}^{\pm}(S)=-{\rm Spec}^{\pm}(S)=\sqrt{-1}\ {\rm Spec}^{\mp}(S)$;
\item ${\rm Spec}^+(S)\subset \sqrt{-1}\ \mathbb{R}$ and ${\rm Spec}^-(S)\subset\mathbb{R}$;
\item If $v\not\in\mathcal{N}$, then $S(v)$ is Jordan simple, i.e. the minimal polynomial of $S(v)$ decomposes as a product of mutually  different irreducible factors.
\end{enumerate}
\end{lemma}
 
Let $S\in\mathcal{P}_n$ for some $n\ge 1$. We set
$${\rm Spec}^{-}(S):=\{0,\lambda_1,-\lambda_1,\ldots,\lambda_l,-\lambda_l\}, $$
where $\lambda_i>0$. It follows that the polynomial
$$\mu_-(X)=X(X^2-\lambda_1^2)\ldots(X^2-\lambda_l^2)$$
is the minimal polynomial of the operator $S(v)$ for all unit timelike vectors $v$. Likewise, the polynomial $\mu_+(X)=X(X^2+\lambda_1^2)\ldots(X^2+\lambda_l^2)$ is the minimal polynomial of $S(v)$ for all unit spacelike vectors $v$. Let $\sigma_k(\lambda)$ denote the $k$-th elementary symmetric function evaluated at $\lambda_1^2,\ldots,\lambda_l^2$. We see from the homogeneity of $S$, i.e. the property $S(Kv)=K^{2n+1}S(v)$, that  
the identity 
\begin{equation}
\label{eqny}
S(v)^{2l+1}+\ldots+\sigma_k(\lambda)(v,v)^{(2n+1)k}S(v)^{2l-2k+1}+\ldots+\sigma_l(\lambda)(v,v)^{(2n+1)l}S(v)=0
\end{equation}
holds for all  $v\not\in\mathcal{N}$.
In fact, it holds for all $v\in T_PM$ by continuity.
\medbreak
We now assume $l\ne0$, i.e. $S\not\equiv 0$, and consider the operator 
$$A(v):=S(v)^{2l}+\ldots+\sigma_k(\lambda)(v,v)^{(2n+1)k}S(v)^{2l-2k}+\ldots+\sigma_l(\lambda)(v,v)^{(2n+1)l}Id.$$
It follows from relation (\ref{eqny}) that $\mathrm{Im}\big{(}S(v)\big{)}\subset\mathrm{Ker}\big{(}A(v)\big{)}$. 
In fact, most of the time we have $\mathrm{Im}\big{(}S(v)\big{)}=\mathrm{Ker}\big{(}A(v)\big{)}$, as shown in the following Lemma.

\begin{lemma}
Adopt the notation established above and let $v\not\in\mathcal{N}$. Then $\mathrm{Im}\big{(}S(v)\big{)}=\mathrm{Ker}\big{(}A(v)\big{)}$.
\end{lemma}

\medbreak\noindent{\bf Proof.}
Let $x\in\mathrm{Ker}\big{(}A(v)\big{)}$. Since the operator $S^2(v)$ is diagonalizable, we have $$V=\mathrm{Im}\big{(}S^2(v)\big{)}\oplus\mathrm{Ker}\big{(}S^2(v)\big{)}.$$
Since $S(v)$ is Jordan simple, we have $\mathrm{Ker}\big{(}S^2(v)\big{)}=\mathrm{Ker}\big{(}S(v)\big{)}$.
Hence $V=\mathrm{Im}\big{(}S(v)\big{)}\oplus\mathrm{Ker}\big{(}S(v)\big{)}.$

Now write $x=x_1+x_2$, where $x_1\in \mathrm{Ker}\big{(}S(v)\big{)}$ and $x_2\in \mathrm{Im}\big{(}S(v)\big{)}$; notice that the inclusion $\mathrm{Im}\big{(}S(v)\big{)}\subset\mathrm{Ker}\big{(}A(v)\big{)}$  implies $x_1=x-x_2\in\mathrm{Ker}\big{(}A(v)\big{)}$. We see from $x_1\in\mathrm{Ker}\big{(}S(v)\big{)}$ that 
$$A(v)x_1=\sigma_l(\lambda)(v,v)^{(2n+1)l}x_1.$$
Consequently, $x_1=0$ and $x=x_2\in\mathrm{Im}\big{(}S(v)\big{)}$.
\qed
\medbreak
Rewriting $\mathrm{Im}\big{(}S(v)\big{)}$ as $\mathrm{Ker}\big{(}A(v)\big{)}$ is beneficial for the following reason. Suppose there exists a vector bundle $E$ over $T_PM-\{0\}$ such that $E\big{|}_v=\mathrm{Im}\big{(}S(v)\big{)}$ for all $v\not\in\mathcal{N}$. Let then $x$ be a section of $E$ over an open set $\mathcal{B}$. Since $\mathrm{Im}\big{(}S(v)\big{)}=\mathrm{Ker}\big{(}A(v)\big{)}$ for all $v\in\mathcal{B}-\mathcal{N}$, we see that $A(v)x(v)=0$ for all $v\in\mathcal{B}-\mathcal{N}$. In fact, it follows by continuity that 
$A(v)x(v)=0\mathrm{\ \ for\ all\ \ } v\in\mathcal{B}.$

It is for this reason that we study the set 
$$\Delta=\Big{\{}x\in\mathbb{R}[T_PM,T_PM]\Big{|} A(v)x(v)=0 \mathrm{\ for\ all\ } v\in T_PM\Big{\}}.$$
Roughly speaking, the elements of $\Delta$ represent the sections of the desired vector bundle $E$. Observe that $\Delta$ is a $\mathbb{R}[T_PM,\mathbb{R}]$-submodule of $\mathbb{R}[T_PM,T_PM]$. The submodule $\Delta$ has the following two properties. 

\begin{lemma}
\label{Delta}
Adopt the notation established above. 
\begin{enumerate}
\item Let $x\in\mathbb{R}[T_PM,T_PM]$ be such that for some non-zero element $P\in\mathbb{R}[T_PM,\mathbb{R}]$ we have $P\cdot x\in\Delta$. Then $x\in\Delta$.
\item Let $T:T_PM\to T_PM$ be the differential of a local isometry  of $M$ preserving $P$. Assume $x\in\Delta$. Then $Tx\in\mathbb{R}[T_PM,T_PM]$ defined by $$(Tx)(v):=Tx(T^{-1}v)$$
is also an element of $\Delta$.
\end{enumerate}
\end{lemma}

\medbreak\noindent{\bf Proof.}
Assume $P\cdot x\in\Delta$ for some non-zero $P\in\mathbb{R}[T_PM,\mathbb{R}]$ and $x\in\mathbb{R}[T_PM,T_PM]$. Then
$$0=A(v)P(v)x(v)=P(v)\cdot A(v)x(v) \mathrm{\ \ for\ all\ \ }v\in T_PM.$$
It follows that $A(v)x(v)=0$ for all $v$ such that $P(v)\ne0$. By continuity $A(v)x(v)=0$ for all $v\in T_PM$. In other words $x\in\Delta$.

Now let $T$ be the differential of a local isometry  of $M$ preserving $P$. Since $S\in\mathcal{P}$, we have $S(Tv)=T\circ S(v)\circ T^{-1}$ for all $v\in T_PM$. As a consequence, $A(Tv)=T\circ A(v)\circ T^{-1}$ and
$$A(v)(Tx)(v)=A(v) Tx(T^{-1}v)=T A(T^{-1}v)x(T^{-1}v)=T(0)=0,$$ 
proving that $Tx\in\Delta$.
\qed
\medbreak
To get our vector bundle we consider the set 
$E\big{|}_v:=\big{\{}x(v)\big{|}x\in\Delta\big{\}},$
i.e. the image of $\Delta$ under the evaluation map $x\mapsto x(v)$. As such, $E\big{|}_v$ is a subspace of $T_PM$ for all $v\in T_PM$. We now study the rank of $E\big{|}_v$. 

\begin{lemma}
Adopt the notation established above. The rank of $E\big{|}_v$ is equal to $r$ for all non-zero vectors $v\in T_PM$.
\end{lemma}

\medbreak\noindent{\bf Proof.}
Note that if $T$ is the differential of a local isometry  of $M$ preserving $P$, then 
\begin{equation}
\label{eqrank}
E\big{|}_{Tv}=T\big{(}E\big{|}_v\big{)}.
\end{equation}
Indeed, let $a\in E\big{|}_v$. This means that $a=x(v)$ for some $x\in\Delta$. Then 
$$Ta=Tx(v)=Tx(T^{-1}Tv)=(Tx)(Tv).$$
By the previous Lemma $Tx\in\Delta$ and so $Ta\in E\big{|}_{Tv}$. We now have 
$T\big{(}E\big{|}_v\big{)}\subset E\big{|}_{Tv}$ as well as $T^{-1}\big{(}E\big{|}_{Tv}\big{)}\subset E\big{|}_v$. Consequently, $T\big{(}E\big{|}_v\big{)}=E\big{|}_{Tv}.$ 
Since our manifold is locally isotropic, the relation (\ref{eqrank}) implies that the rank of $E\big{|}_v$ is independent of the choice of unit spacelike (resp. unit timelike, non-zero null) vector $v$.

In fact, we know that the rank of $E\big{|}_v$ for unit spacelike and unit timelike vectors $v$ is exactly $r$. This follows from 
\begin{equation}
\label{eqrank2}
E\big{|}_v=\mathrm{Im}\big{(}S(v)\big{)} \mathrm{\ \ for\ all\ \ } v\not \in\mathcal{N}.
\end{equation}
To prove this equality we first consider $a\in E\big{|}_v$. Since $a=x(v)$ for some $x\in\Delta$, we have $A(v)a=A(v)x(v)=0$, i.e. $a\in\mathrm{Ker}\big{(}A(v)\big{)}=\mathrm{Im}\big{(}S(v)\big{)}$. Conversely, take $b\in\mathrm{Im}\big{(}S(v)\big{)}$. We have $b=S(v)b^\prime$ for some $b^\prime\in T_PM$. The map
$$x:w\mapsto S(w)b^\prime$$
is an element of $\Delta$ with $x(v)=b$. In other words, $b\in E\big{|}_v$.

Let $x_1,\ldots,x_r\in\Delta$ be such that the vectors $x_1(v),\ldots, x_r(v)$ are linearly independent at some $v\in T_PM-\mathcal{N}$. Using Lemmas \ref{ind1}, \ref{ind3} and  \ref{Delta} we see that $I(x_1,\ldots, x_r)\ne\{0\}$ and hence there exist $y_1,\ldots, y_r\in\Delta$ with  $I(y_1,\ldots,y_r)\not\subset\mathcal{I}$.  Lemma \ref{ind2} implies that $y_1(v),\ldots,y_r(v)$ are linearly independent for some $v\in\mathcal{N}$ and therefore the rank of $E\big{|}_v$ for $v\in\mathcal{N}$ is at least $r$.

On the other hand, the rank of $E\big{|}_v$ cannot be bigger than $r$. This follows from the fact that if $y_1(v),\ldots,y_{r+1}(v)$ are linearly independent at $v_0\in\mathcal{N}$, they are linearly independent on a neighborhood $\mathcal{B}$ of $v_0$. Since the rank of $E\big{|}_v$ is $r$ for all $v\in\mathcal{B}-\mathcal{N}\ne\emptyset$, the vectors $y_1(v),\ldots,y_{r+1}(v)$ cannot be linearly independent for all $v\in\mathcal{B}$.

Consequently, the rank of $E\big{|}_v$ is $r$ for all $v\in T_PM-\{0\}$.
\qed
\bigbreak
It now follows that the disjoint union $$E:=\bigcup E\big{|}_v$$ is a vector bundle over $T_PM-\{0\}\simeq S^{p+q-1}$ of rank $r$. Local triviality of $E$ around $v_0\in T_PM-\{0\}$ is ensured by polynomial maps $x_1,\ldots,x_r$ such that $x_1(v_0),\ldots,x_r(v_0)$ are linearly independent. We see from relation (\ref{eqrank2}) that the vector bundle $E$ has the desired property $E\big{|}_v=\mathrm{Im}\big{(}S(v)\big{)} \mathrm{\ for\ all\ } v\not \in\mathcal{N}$.

Like most of the vector bundles we have encountered so far, the vector bundle $E$ descends to a vector bundle over the real projective space. This is due to the following observation.

\begin{lemma}
Adopt the notation established above. We have $E\big{|}_v=E\big{|}_{-v}$ for all $v\in T_PM-\{0\}$.
\end{lemma}

\medbreak\noindent{\bf Proof.}
Let $a\in E\big{|}_v$. Then $a=x(v)$ for some $x\in\Delta$. Now consider $y\in\mathbb{R}[T_PM,T_PM]$ defined by $y(w):=x(-w)$. Relation $A(w)=A(-w)$ implies  
$$A(w)y(w)=A(-w)x(-w)=0, \mathrm{\ \ i.e. \ \ } y\in\Delta.$$
Consequently $a=y(-v)\in E\big{|}_{-v}.$ We now see that $E\big{|}_v\subset E\big{|}_{-v}\subset E\big{|}_v, \mathrm{\ \ i.e.\ \ }E\big{|}_v=E\big{|}_{-v}.\qed$
\medbreak
We shall use $\mathrm{Im}(S)$ to denote the vector bundle over $\mathbb{R}P^{p+q-1}$ induced by $E$. Note that if $\iota_2:\mathbb{R}P^{q-1}\to\mathbb{R}P^{p+q-1}$ is the natural inclusion induced by the inclusion of $\mathbb{R}^q$ as a maximal positive definite subspace of $T_PM$, then 
$$\iota_2^*\big{(}\mathrm{Im}(S)\big{)}\cong\mathrm{Im}(S)_2.$$
In particular, we have $\iota_2^*\big{(}\mathrm{Im}(S)\big{)}\cong\iota_2^*\big{(}\mathrm{Im}(S)\big{)}\otimes\gamma_1$.
\bigbreak
We are now ready to prove Wolf's Theorem.

\begin{theorem}
Let $M$ be a locally isotropic pseudo-Riemannian manifold of signature $(p,q)$, where $p\ne q$ and $\max\{p, q\}\ge 11$. Then $M$ is locally symmetric.
\end{theorem}
\medbreak\noindent{\bf Proof.}
Adopt the notation established throughout this Section.  As pointed out earlier, it is enough to consider the case when $p\le q-1$. Let $\mathbf{S}$ be the Szab\'o operator at $P\in M$. By Theorem \ref{thm2} it is enough to show $\mathbf{S}\equiv 0$. We assume the opposite and argue for contradiction.

If $\mathbf{S}\not\equiv 0$, we may perform the construction explained above. Therefore, if $\mathbf{S}\not\equiv 0$, we have a vector bundle $\mathrm{Im}(\mathbf{S})$ over $\mathbb{R}P^{p+q-1}$ with non-zero rank and such that:
\begin{itemize}
\item $\mathrm{Im}(\mathbf{S})$ is a sub-bundle of the trivial vector bundle $\mathbb{R}P^{p+q-1}\times T_PM$;
\item  $\iota_2^*\big{(}\mathrm{Im}(\mathbf{S})\big{)}\cong\iota_2^*\big{(}\mathrm{Im}(\mathbf{S})\big{)}\otimes\gamma_1$, where $\iota_2:\mathbb{R}P^{q-1}\to\mathbb{R}P^{p+q-1}$ is the natural inclusion. 
\end{itemize}
Since by assumption $q-1\ge 10$ and $\frac{p+q-1}{2}\le q-1$, Lemma \ref{techn} implies the rank of $\mathrm{Im}(\mathbf{S})$ is $0$. Contradiction.
\qed

\section{Final Remarks and Acknowledgments}

\hskip .2inLet $M$ be a locally isotropic manifold of signature $(p,q)$, where at least one of the integers $p$ and $q$ is odd. It is known that the manifold $M$ has constant sectional curvature, see \cite{Wolf}. As pointed out in the Introduction, there is an elementary proof of this constant sectional curvature result in the Lorentzian signature. The proof uses the Jacobi operator and the methods related to the ones we use in this paper. It seems likely that the application of the polynomial methods developed in Sections 3 and 5 will lead to new proofs of the constant sectional curvature result. It is also possible that   "gluing over the nullcone" which we performed in Section 5 is a special case of a more general phenomenon having further consequences. 

This paper is dedicated to Peter, Paul, Peter, Alan and Amelie.

\medbreak\noindent{\bf Address:} Iva Stavrov, 

Department of Mathematical Sciences, Lewis and Clark College, 

0615 SW Palatine Hill Road, Portland, OR 97219, USA

\noindent{\bf email:} istavrov@lclark.edu
 \end{document}